\documentclass[11pt,twoside]{atmp} % \documentstyle[amscd,amssymb,verbatim,12pt]{amsart}

\usepackage{amsmath,amssymb,latexsym}

\input xypic \xyoption{curve} %\xyoption{all}
\usepackage{color}

\copyrightnotice{2007}{1}{1}{22} %% year, volume, first page, last page.

\newcommand{\C}[1]{{\mathcal#1}} % Calligraphic
\newcommand{\F}[1]{{\frak#1}}% Fraktur

\newcommand{\U}{{\C{U}}}% 

\begin{document}

\title[The Feynman Legacy]{The Feynman Legacy}
\arxurl{TBA} % <hep_reference_#>}

\author[L. M. Ionescu]{Lucian M. Ionescu}

\address{Illinois State University, Mathematics Department, IL 61790-4520}  %lines should be separated with double backslashes: \\
\addressemail{LMIones@ilstu.edu}

\begin{abstract}
The article is an overview of the role of 
graph complexes in the Feynman path integral quantization.

The underlying mathematical language is that of PROPs and operads,
and their representations.

The sum over histories approach, the {\em Feynman Legacy},
is the bridge between quantum physics and quantum computing,
pointing towards a deeper understanding of the fundamental concepts
of space, time and information.
\end{abstract}

\maketitle

%=================================================================
%                   Introduction
%=================================================================
\section{Introduction}\label{S:intro}
We present the role of graph complexes in deformation quantization,
focusing on the work of M. Kontsevich on deformation of Poisson manifolds \cite{K-dq}
and the author's joint work with D. Fiorenza \cite{FI}.
Their role for the Feynman Path Integral quantization method  is emphasized \cite{CF}.

The Feynman Legacy, a natural amplification of the original formulation 
of Quantum Mechanics (QM) by W. Heisenberg, 
is presented in the context of the various 
quantum theories (QFT, CFT, string theory etc.), as an interplay 
between the categorical structures (operads, PROPs etc.) and
the computer science interpretation: 
automata as ``states and transitions'' approach to modeling.

The article is based on the recent talks  on the subject \cite{I-talks}.
It is not a comprehensive overview, 
rather aiming at a few important aspects, pertaining to the interpretation
and motivation of the present theories, and therefore subjective.

\cutpage %move this line so that the first page breaks at the appropriate place.

\setcounter{page}{2} % <insert page # for second page>}

\noindent

\vspace{.2in}
{\bf Acknowledgments.}
The article is based on an overview talk given at Bucharest University,
and a second talk at University of Rome, La Sapienza.

I would like to thank Gabriel Pripoae, the organizer of the Gh. Vranceanu seminar
at Bucharest University, and Domenico Fiorenza
for the invitation to speak in the Algebraic Geometry seminar.

I would also like to thank for the impressive research and lodging conditions
while visiting I.H.E.S., where some of the present ideas matured and contributed
to a better ``big picture'' of the subject matter in the second talk.

% ***********************************************************************
%			Past: Deformation Quantization
% ***********************************************************************
\section{Past: deformation quantization}
The initial formulation of QM by W. Heisenberg was specifically 
formulated as an automaton model,
in terms of ``states and transition amplitudes'', 
which can be viewed as a complexified version of Markov processes
(yet with a continuum time).

The holistic, intrinsic formulation in terms of Hilbert spaces and unitary transformations,
of von Neumann tended to blur this ``in coordinates formulation''.

Later on, to tame QM, i.e. to make it as classical as possible,
M. Flato, A. Lichnerovich a.a. have launched the program of {\em deformation quantization},
as a ``minimal'' change to the classical formulation of mechanics,
to make it ``quantum'': start from the Poisson algebra of 
{\em classical observables} on a Poisson manifold, 
and ``just'' make it non-commutative to accommodate 
Heisenberg's uncertainty (commutation) relations!
In our opinion, in this way the ``spirit'' of QM was frozen ... 

On the other hand, in QFT, a new revolutionary approach to quantum physics
was introduced by Feynman (following a suggestion of Dirac!).
It was well established amongst physicists, 
being now considered the most powerful quantization method,
yet only hardly tolerated by mathematicians as being not ``rigorous''
\footnote{Not ``digitally signed'' by mathematicians.}

Yet its power was reveled one more time by Kontsevich's solution to the
deformation quantization problem, as explained in \cite{CF}.

\subsection{How to build a star product?}
Kontsevich's formula
$$f\star g=\sum_\Gamma W_\Gamma B_\Gamma,$$
at ``top level'' looks quite natural.
It is a sum over graphs of bidifferential operators $B_\Gamma$ associated
to graphs $\Gamma$ according to a rule similar to a 
Feynman rule from QFT (perturbative approach).
When more details are revealed, one starts to wander why it really works!
Given a Poisson manifold $(M,\pi)$, with its Poisson algebra of observables
$A=C^\infty(M)$, vector fields $\F{g}=Der(A)$, polyvector fields
$T_{poly}=\bigwedge^\bullet\F{g}$ and differential operators
$D_{poly}=Hoch^\bullet(A;A)$, 
Kontsevich defines a pairing between graphs and polyvector fields
which in particular yields the above bidifferential operators.
The above product obviously perturbed the multiplication,
since $B_{b_0}$ is by definition the commutative multiplication
and $B_{b_1}=\pi$ is the Poisson bracket 
(see \cite{DF} for notations; a detailed account appears in \cite{CKTB} etc.).
The ``problem'' was to define the coefficients $W_\Gamma$ to ensure 
associativity;
they are defined by Kontsevich in a similar (cryptic) manner, 
using the same graphs,
and a Feynman rule with a closed 2-form $\alpha$ derived from the
angle form on the Poincare half-plane.
Amazingly, the formula works (!), and
provides not only a star-product, but also an $L_\infty$-morphism
proving the Formality Conjecture ... but again, why?

\subsection{Why the formula ``works''}
On the physics side, their result (deformation quantization), 
is the outcome of a deeper theory: Feynman Path Integral quantization method;
physicists know it is the most powerful quantization method available.
Form such a point of view, the deformation quantization goal is guaranteed:
quantization leads to quantization!

The string theory (sigma model) interpretation of Kontsevich's construction
was reviled in \cite{CF}.

From the mathematics side, this result (formality) hides 
a deeper algebraic structure: 
the DG-coalgebra (Hopf algebra) of graphs  and cocycle condition,
or from a dual point of view the DG-Lie algebra of graphs and 
Maurer-Cartan (MC) equation of deformation theory, as disclosed
by the present author in \cite{I-pqft,I-cfg,I-comb}.
The result was later explained in a ``round form'' in \cite{FI}:

1) The Kontsevich rule $\C{U}$ is a morphism of DGLAs,

2) the coefficients $W$ yield a DG-coalgebra 1-cocycle 
(MC-solution of the dual DGLA),

3) $\C{F}=W\C{U}$ is an $L_\infty$-morphism.

The rest is a standard argument:

4) The Poisson structure $\alpha$ being a MC solution in the
Lie algebra of polyvector fields (Jacobi identity), 
is mapped by the $L_\infty$-morphism to a MC-solution,
providing by exponentiation a perturbation of the commutative product:
$$\star=\mu+\C{F}(e^\alpha).$$
For the present author the ``underlying ideas'' (algebraic structures involved)
are more interesting then the result itself
\footnote{I is a pity ``Part II'' was never written!},
especially in connection with renormalization in the framework of 
Connes and Kreimer (see \cite{IM}).

%=================================================================
%                   Present: ren and graph homology
%=================================================================
\section{Present: renormalization and graph homology}
During the author's first ``IHES period'',
working at the same time on renormalization and
deformation quantization turned out beneficial,
leading to the ``unification'' of Kremer's coproduct and
Kontsevich's homology differential.
With hindsight, glancing at Kontsevich's proof it is hard to miss
that the sums of products must be a coproduct and that the one edge case
belongs to the homology differential.
``Reverse engineering'' Kontsevich's result
\footnote{... concepts, concepts, concepts!}
and pulling everything back to the source, the graphs,
yields the definition of the coproduct $\Delta$ and that
$W$ satisfies the 1-cocycle equation $W(\Delta\Gamma)=0$.
The rest are technicalities: 1) reformulate Kreimer's coproduct
and applying it to Kontsevich's formula, 2) include Kontsevich's 
differential $d$ as a part of the coproduct \cite{I-pqft}.

\subsection{The homological algebra interpretation}
At this point we see the shadow of the cobar construction $D=d+\Delta$!
Now, turning ``on'' the homological algebra machinery
yields graph cohomology \cite{I-cfg}.

Besides the nice mathematical result, a deeper idea emerges:
the ``resolution of degrees of freedom'' is a crucial mathematical idea
standing for the physicist's {\em multi-scale analysis},
as it will be explained later on.

\subsection{Main results explained}
The first main question to be answered in order to solve ``the Kontsevich's puzzle'' 
was ``What is the map $\U$, really?''.
In \cite{FI} we showed that it is a DGLA morphism,
and therefore mapping the Maurer-Cartan solutions 
in the DGLA of graphs to the Eilenberg-Chevalley (EC) complex,
i. e. $\C{U}$ maps the MC-solution $W$ to a MC-solution
$$U=\C{U}(\sum_\Gamma W_\Gamma\Gamma).$$
Such an MC-solution is an L-infinity morphism
and also a quasi-isomorphism (since the degree 0 part is),
therefore concluding Kontsevich's proof of the {\em Formality Conjecture}.

The conceptual (mathematical) interpretation of Kontsevich's solution
of the formality conjecture and deformation quantization
of Poisson manifolds \cite{FI} includes the following results:
1) The Kontsevich graphs $\C{G}$ have a {\em DG-coalgebra} structure,
or equivalently a (dual) DGLA structure (\S 5);
2) The Kontsevich graphical calculus
$\C{U}$ is a DGLA morphism (\S 6) ({\em calculus of derivations});
3) If $W$ is a cocycle of $\C{G}$ (\S 7),
$\C{U}(W)$ is an L-infinity morphism ({\em the formality morphism}) (\S 4);
4) $\C{F}(exp(\pi))$ is a {\em star product} (\S 4);
5) There is a solution involving only 
{\bf graphs without circuits} (\S 8), 
we call a {\em semi-classical star product}.
A few details will be given next,
to later justify the framework we chose for what we call Feynman Processes.

The {\em DG-coalgebra of Feynman graphs} (\S 5)
$\C{G}^{\bullet,\bullet}$ is a bigraded DGLA,
associated to a {\bf pre-Lie algebra of graphs}
\footnote{The pre-Lie structure is fundamental concept;
it is a {\em formal connection},
due to the properties it enjoys \cite{I-tor}},
or dually a coalgebra, with coboundary differentials.
The {\bf internal differential} is Kontsevich's homology differential,
and the {\bf external differential} is Hochschild differential 
corresponding to Gerstenhaber composition by a prefered element.
Such a ``bi-DGLA'' will be the object of study of a different article,
focusing on the resolution of a complex, or in categorical language,
the resolution of a PROP.

The Kontsevich rule defining ``a la Feynman'' $\U$,
the {\em DGLA morphism $\C{U}$} (\S 6),
is a {\bf graphical calculus of derivations} (\S 2):
vector fields $\partial_i$ act on functions $f$:
$\overset{X}{\to}\overset{f}{\bullet}$.

{\em The formality morphism} between the two DGLAs: 
$$\C{F}=W\C{U}:T_{poly}\to D_{poly}$$
is an L-infinity morphism, 
because $W$ is a cocycle in the DG-coalgebra of graphs
(solution of the MC of the dual DGLA of graphs).

Since its zero degree is known to be a quasi-isomorphism
(Kostant-Hochschild-Rosenberg Theorem),
it proves that the $D_{poly}$ DGLA is {\em formal},
i.e. quasi-isomorphic to its cohomology.

A few more words about the DG-coalgebra cocycle $W$ are in order.
The cocycle refers to the cobar construction associated to the
DG-coalgebra \cite{I-cfg}:
$$0=(D^*_{cobar} W)(\Gamma)=W(d\Gamma+\Delta_b\Gamma)$$
$$=\sum_{e\ internal\ edge}\pm W_{\Gamma/e}+\sum_{\gamma\subset\Gamma\ boundary}
\pm W_\gamma W_{\Gamma/\gamma},$$
It can be rewritten as the {\bf Maurer-Cartan solution} in the dual DGLA:
$$=\delta W+\frac12[W,W](\Gamma).$$
It is defined using a similar ``Feynman rule''
with (winding number/angle form) $\alpha=d\theta$ defining 
the bivector field $dA=d\theta\wedge d\theta$.

The {\em star product} is obtained as a deformation of the 
commutative multiplication by a perturbation.
To ensure associativity, the perturbation
must be a MC-solution $\C{F}$ in the CE-complex:
$$Hom(T^\bullet(T_{poly},D_{poly})\cong Hom(T^\bullet(\F{g},D_{poly}).$$
Now $\C{F}=W\C{U}$, 
maps the exponential of a Poisson structure,
i.e. a  MC-solution in $T_{poly}$ (with trivial differential):
$$d\pi+[\pi,\pi]=[\pi,\pi]=0\quad (Jacobi\ identity),$$
to a MC-solution in $D_{poly}$, yielding the desired 
{\bf perturbation} $\partial$, of the commutative product:
$$\star=\mu+\partial,\quad \partial=\C{F}(exp(\pi)).$$
Conform the standard philosophy of deformation theory,
MC-solutions modulo equivalence,
is thought of as the tangent space at $\mu$ at the
the moduli space of deformations of associative algebras 
with given underlying vector space.

Two questions were naturally asked at this point:
1) Are all graphs necessary?, and 
2) Is there a simpler solution $W$?

If considering graphs without circuits, i.e. forests $\C{F}$,
the results hold essentially because the inclusion
$\C{F}\subset \C{G}$ is a quasi-isomorphism of complexes,
yet a different cocycle $W$ must be provided.
The fact that in Kontsevich solution the angle form
may be replaced with any closed form,
is an indication that the coefficients may be universal in some sense.
Indeed in \cite{I-comb} it is argued that such a simple/universal solution
exists, reminiscent of the Hausdorff series, universal with respect
to the particular Lie algebra considered.
A detailed proof is currently under consideration \cite{IS}.

% *******************************************************
% 		Future
% *******************************************************
\section{Future: Quantum Information Dynamics}
Before being more specific, a brief motivation is in order.

\subsection{Feynman Processes and quantum computing} 
The {\bf method} of {\em Feynman Path Integrals}, 
is a way of ``thinking''; it is modeling in terms of
{\em states} and {\em transitions}, i.e. {\em automata}
in the sense of mathematics and computer science.
It can be applied to {\em classical mechanics} as well as to QM/QFT.
(recall that QM may be presented as a (0+1)-QFT).

We would like to capture this conceptual bridge
as the ``motto'': 
$${\bf Quantum\ Physics}\ IS\ {\bf Quantum\ Computing}$$
a point of view in fact adopted by Feynman himself.

Expanding the idea:
$${\bf Quantum\ interactions\ are\ Quantum\ Communications}$$
it enables a unified Mathematical-Physics and Computer Science description,
without the asymmetry {\em system $\to$ observer},
towards the incorporation of {\bf entropy} and {\bf information}
as part of the foundations of physics.
Physicists know that quantum information is more fundamental
than matter (fields included),
for example at the level of the {\em Information Paradox}
in the context of a radiating black hole \cite{NS-ip}.

An extended {\em Equivalence Principle} between energy and information
was conjectured as mandatory to really progress beyond the
present attempts to unifying gravity and QM,
not just as quantum gravity, but including all known fundamental
interactions \cite{I-UP}.

As an example of line of thought in this direction, 
the ``unit of a black-hole surface'' is thought of as a bit of space-time 
by Lee Smolin \cite{S-QG}; 
more that, we claim, 
quantum black holes should be considered ``generic objects'',
towards a spin/qubit network (spin foams \cite{R1,R2} etc.) 
description of quantum physics.

\subsection{Feynman Processes: the mathematical-physics}
The above ``physics interface'' has already an implementation
underway. 
It iterates the ``old'' idea of representing
a non-linear structure,
i. e. a functorial representation of a tensor category,
as a replacement of the mechanics point of view centered on the 
concept of {\em manifold},
where the linearization is obtained at the level of the tangent space,
therefore binding rigidly the configuration/state space and 
the symmetries/dynamics.

A {\em representations} of a {\bf causal structure}
{\em geometric category} playing the role of ``space-time'',
i. e. taking into account the external degrees of freedom,
in some sense the macroscopical structure of the system 
(NO embeddings in an ambient space yet!).
As examples, we have the (PROP or Hopf algebra of) 
Feynman-Kontsevich graphs, Segal PROP, 
cobordism categories, lattices etc.
The causal structure captures the data regarding ``what influences what'',
and the norm is that there is no nice cartesian product of space
versus time, not even locally, especially since one has to abandon
locality to really account for quantum phenomena.

Then a {\bf functorial representation} is just an algebra over the PROP:
QFT, CFT, ST, TQFT etc.
We will later endow such a PROP with additional structure,
and call it a {\em Feynman Category}.

\subsection{Physics interpretation and motivation}\label{S:physmotiv}
As mentioned above,
in a {\bf causal structure} there may be no locally defined space-time structure!
(There may be loops - ``quantum feedback'' etc.).

The causal structure is {\bf NOT ``fixed''}, like a manifold is, 
but it is a {\bf resolution},
allowing to implement the multi-scale analysis, 
either in a homological algebra fashion or similar to the 
multi-resolution analysis of Haar wavelets etc.
the {\bf processes} are inter-related as a {\bf complex};

From this point of view,
the approach is NOT a ``{\bf perturbative approach}'',
but it may be obtained as the {\em outcome} of one:
$$\diagram 
Gaussian\ integral \rrto & & Feynman\ graphs\ and\ rules\\
Z=\int\C{D}\phi e^{S(\phi)}\quad \rrto^{Wick's Th.} \quad 
& &  \quad \sum_\Gamma \C{F}(\Gamma)/|Aut(\Gamma)| h^{deg(\Gamma)} \
\enddiagram$$
It is similar to the dual role and origin of the exponential:
$$e^x=\sum_n x^n/n!:\quad \frac{dy}{dx}=y\quad or\quad Aut([n])?$$
as an analytic solution of a smooth differential equation,
or as a generating function of the automorphism groups of 
finite sets (categorification);
one can obtain the categorification perturbatively
as a Taylor expansion ...

Our point of view is that 
FPI and {\bf matrix-models} {\em are}  ``{\bf recipes}'' for 
representations of causal structures.

No matter the ``origin'' or interpretation,
the productive language for modeling quantum phenomena is that of
PROPs and operads.

% ***********************************************************
%			From groups to Operads and PROPs
% ***********************************************************
\section{From groups to operads and PROPs}
PROPs and operads are enriched tensor categories.
Familiar examples for physicists are Feynman graphs, 
Riemann surfaces, and even lattices or almost any other
collection of discrete models, can be presented using this language,
suited for an ``automaton/computer science'' interpretation.

This specialized categorical language has a big advantage over 
other presentations: 
it is a ``graphical interface'' to Quantum Mathematical-Physics.

\subsection{PROPs and operads}
If we look back, 
abstract groups were defined long after the use of groups of
transformations in mathematics and physics,
i. e. after the representation theory started to develop.

Groups may be thought of as a one object category,
encoding the algebraic structure of 
groups of transformations,
i.e. the composition law with its basic properties:
associativity, inverse etc.

Operads, in a similar way,
encode various types of algebraic structures,
like the classical well known ``binary operation spaces'' (monoids, algebras etc.)
of various types: associative, commutative and Lie.

A short formal definition is the following (see \cite{DF} for details): 
a $\C{C}$ operad over $N$ consists of 
``Input/Output'' operations $\C{O}(n)\in Ob(\C{C})$ 
and rules for composing these operations (the motherboard model; see \cite{Leisner}).

Then the usual theory can be developed:
ideals representing ``constraints'' and quotients by these ideals
yielding new examples implementing the constraints.

We then have presentations by {\em generators and relations},
as in the following examples.
We will ignore the important aspect of the symmetric action,
mentioning the non-sigma version of operad/PROP (PRO - see \cite{MSS}).

\subsection{Examples}
\subsubsection{The operad {\bf Assoc} }
The free operad generated by the graphical symbol of a generic
binary operation is the operad of binary trees.

The quotient of this operad by the ideal generated by one element
representing in a graphic/schematic way the associativity constraint
is again an operad, the {\em Associative Operad}.
What is then an associative algebra, in this language?
To transform the ``generic'' into ``actual'',
we need instantiation: give meaning to the graphical symbol
of the generic binary composition.
In various other areas this operation is called {\em coloring}
the graphical symbol with ``colors'',
which do not have to be in fact colors, nor even numbers;
the colors (labels) may be elements of objects of categories,
and usually are in fact set theoretic operations.

This coloring process is a morphism of operads,
therefore qualifying for the term of representation,
since the target enriched category is usually k-linear.
So, an associative algebra is a {\em representation}
of the Associative Operad.
Because of the above tight connection,
a representation of an operad $\C{O}$ is also called 
an {\em algebra over the operad $\C{O}$}.

\subsubsection{The operad {\bf Lie}}
The other classical examples are quite similar,
differing only at the level of the generator of
the quotient ideal of the free operad of binary trees.

For the {\em Lie Operad},
take a graphical (operational) representation of the
Jacobi condition as one generator of the ideal of the
free operad of binary trees.
The other generator represents the anti-commutativity
(here we need sigma-operads, with an action of the
symmetric group defined on the objects underlying 
the tensor category of an operad).

Lie algebras are now algebras over the Lie Operad.

\subsubsection{The operad {\bf Comm}}
To realize commutative algebras as algebras
over the {\em Commutative Operad},
quotient the free binary trees operad by
the ideal generated by the graphical representation
of commutativity and associativity.

\subsection{Representations of operads}
As already exemplified above,
representations of operads,
also called algebras over the corresponding operad,
are morphisms between operads, where the target operad is k-linear:
$\rho:\C{O}\to Vect_k$.
Here the {\em endomorphism operad}
$End(V)(n)=Hom(V^{\otimes n},V)$ is just a colored instance 
of the binary tree operad:
to a node associate a vector space $V$,
to a set of external/terminal nodes (concatenation is the monoidal product),
associate the corresponding tensor product.
To internal nodes associate the possible ``interactions'',
the Hom-set of linear functions, 
and the operad's composition is the usual composition of 
(multi)linear functions.

Now since the range of such a representation of the Associative
Operad is 1-dimensional, generated by the only binary operation
graphical symbol,
it is determined by specifying $V$ and an ``actual'' binary operation,
the binary multiplication on $V$ which must satisfy associativity:
$\cdot:V\otimes V\to V$.

The above example is obviously a graphical reformulation 
of basic linear algebra facts,
in a categorical language with a definite quantum physics flavor,
as one can see comparing with the (1+1)-TQFT case:
with a bit more structure in the geometrical category
(0+1-cobordisms), the data needed to specify a representation
amounts to a {\em Frobenius algebra},
as it will be explained later on.
It is important to realize at this point that all the structure
is abstractly captured in the geometric category
being represented, and this is a stage which must be clarified
before further developing the representation theory.
For example, in the (0+1)-TQFT case,
it was shown in \cite{I-note} that the generator of the geometric
category of 1-dimensional cobordisms {\em is} a 
Frobenius algebra in that monoidal category, if one disregards
k-linearity.
Now, what is a {\em String Theory} without a sigma model nor 
a {\em string background}? (to be continued).

% **********************************************************
% 			PROPS and Operads
% **********************************************************
\section{PROPs and ... beyond!}
The operads are just ``half of the (modeling) story'',
since they are the heritance of classical mathematics
modeling interactions as ``combining inputs'' and then 
trying to solve equations, with an idealistic hope for
uniqueness and determinism; reality is much more interesting than this!
The possibilities unfold as a shower of high energy particles,
in a ``network of instances'' rather closer to fireworks ...
How to model this? ... using {\em coproducts},
which are well known in combinatorics as a marvelous tool to keep track
of the various possibilities to achieve a given configuration.
As a quick example, a coproduct {\em dual} to multiplication $(N,\cdot)$ is 
$$\Delta n=\sum_{d\cdot q=n}(d,q),$$
which is essentially the list of divisors of a natural number 
and the corresponding quotient.

At the categorical level, moving from the classical algebraic holistic structures
to objects and relations oriented framework, prone for an automaton interpretation,
the corresponding dual notion to an operad is a cooperad: just reverse the trees.
This is a gain ``half'' of the picture. The enveloping structure is that
of a PROP (PRoduct Operations and Permutations),
which is an enriched version of a tensor category.
To simplify the exposition, we will disregard the permutation action.

A {\em $\C{C}$-PROP over $\C{N}$}, denote by $P$,
consists of objects of a category $\C{N}$ but with 
spaces of morphisms, {\em Hom}s, from a possibly different category,
with additional structure involved.

Then, there are the usual compatible operations of composing morphisms:
{\em vertically}, the tensoring of objects and morphisms, 
and {\em horizontal composition} of source-target compatible morphisms.

The objects are often in correspondence with natural numbers,
and tensoring corresponds to concatenation,
therefore counting the elementary pieces;
we say it is a PROP over $(N,+)$.

Then the morphisms appear as operations with 
multiple inputs and multiple outputs and compositions: $\C{P}(m,n)$.

Now operads become ``half-PROPs'',
as underlying subclasses of morphisms with only one output: 
$$\C{O}\subset \C{P}, \quad \C{O}(n)=\C{P}(n,1).$$
tensoring and composing these 1-output operations 
yields the {\em standard PROP generated by an operad},
denoted by $\C{P}(\C{O})$.

\subsection{Examples}
The endomorphims PROP $\C{E}nd(V)$ of a k-vector space
consists of the linear mappings:
$$\C{P}(n,m)=Hom_k(V^{\otimes n},V^{\otimes m})$$
together with the usual tensor and composition operations.

Trees form the prototypical operad, 
which generates the forest as the standard PROP.
Imposing various constraints as generators of {\em ideals}
to quotient by the free operad, yields the classical
operads responsible of the various kinds of algebras:
associative, commutative and Lie, as their {\em representations}
(see bellow).

More general classes of graphs can be presented as PROPs.
In QFT the graphs may be restricted to have only three,
or only four etc. vertices; the so called
$\phi^3,\phi^4$ theories.
Such a class of graphs forms also a PROP:
the {\bf Feynman PROP}.

Similarly, the category of cobordisms,
i. e. manifolds with boundary split into a negative part,
{\em the source}, and a positive part,
the target, is a PROP. 
The tensor operation on objects is concatenation,
as in probably all such geometric categories,
and the composition of the morphisms is defined
by prescriptions of gluing cobordisms, usually defined on
representatives and compatible with the pertaining equivalence relations.

If instead of topological surfaces
we take into account a complex structure,
then Riemann surfaces (RS) with boundary (holomorphic disks with orientation,
and identification data),
form a PROP, the {\em Segal PROP}.
The tensor operation is again concatenation of circles 
(boundaries of RS/disks),
while composition is the operation of sewing Riemann surfaces,
i.e. topological gluing and continuation of complex structure.

% ***********************************************************
\subsection{Algebras over PROPs}
Operads and PROPs encode the type of operations had in mind.
The actual examples come as their representations,
which are k-linear valued morphisms of PROPs/operads.
Such a morphism is also called an {\em algebra over the operad/PROP}
in consideration.

A {\em representation}  $\rho:\C{P}\to Vect_k$ of a PROP $\C{P}$,
is a strict tensor functor valued in the category 
of k-linear vector spaces \cite{DF}.
Since the PROP is assumed to be over $(N,+)$,
a tensor category with one generator,
the image is a 1-dimensional PROP, 
generated by a vector space $V$:$\C{E}nd(V)$.
Such a representation amounts to a family of morphisms:
$$\rho_{n,m}:\C{P}(n,m)\to Hom(V^{\otimes n},V^{\otimes m}),$$
compatible with the two compositions.

Note that the endomorphism PROP is also 
an algebra over the free PROP of graphs generated by the $(n,m)-corollas$
(nodes with $n$ inputs and $m$ outputs):
% Question: is it the standard PROP generated by trees?

Various classical examples are presented in \cite{DF}.
The more sophisticated examples of quantum physics 
are discussed in \cite{Vor} % (see A. A. Voronov - hep-th/9401023)

\subsubsection{Quantum Field Theory}
Briefly said, Feynman rules of a given QFT,
may be presented functorially as an algebra over the corresponding
Feynman PROP.

\subsubsection{Topological Quantum Field Theory}
An algebra over the $(n+1)$-cobordism category,
i.e. with boundary manifolds of dimension $n$,
is called a {\em topological quantum field theory} (TQFT) \cite{Atyiah}.

In the special case of $(1+1)$-TQFTs, where the objects are
topological surfaces with a disjoint union of circles as boundary,
the data needed to specify such an algebra, i.e. a (1+1)-TQFT,
is equivalent to a {\em Frobenius algebra} \cite{Kock}.
As noted before,
the {\em abstract structure being represented} 
can be recognized at the level of the geometric category;
it is a {\em Frobenius object}, as explained in \cite{I-note}.

\subsubsection{Conformal Field Theory}
Representing the Segal PROP amounts the constructing a CFT
\cite{Huang}.
The various approaches, whether geometric or analytic,
have in fact underlying algebraic structures of the ``analytic type''.

\subsubsection{String Backgrounds}
A String Theory is built by embedding RS into a given space-time manifold,
called a {\em sigma model}. In contrast to the CFT case,
we refer to this approach as ``smooth type''.

Not to refer to an ambient space-time,
an axiomatic approach may be defined:
a {\em String Background} is a representation of the 
PROP of chains of Riemann Surfaces \cite{Vor}:
$$C_\bullet\C{P}(m,n)\to End(H,Q)(m,n),\qquad Q^2=0.$$
The main point is that a {\em differential structure} is taken into
account.

Yet more structure is needed, which in fact is already present in QFT: 
the Feynman graphs have not only a Hopf algebra structure \cite{Kre-HA},
but also a compatible (homology) differential \cite{I-pqft}.
In the categorical language,
a {\em Feynman Category} is a DG-coalgebra PROP,
modeling not only the states and paths,
but also the multi-scale property of degrees of freedom.
It enables the multi-scale analysis or homological algebra machinery 
(causal resolution \cite{I-DWT}, 
cohomology and underlying quantum theories \cite{Vor} etc.),
encoded as a coalgebra structure on ``paths''.

Then a {\em Feynman Process}, as a generalization of QFT, CFT, ST etc.,
is an algebra over the Feynman category,
as a {\em representation of the causal structure}.

% **********************************************************
% 			Perturbative ... or not?
% **********************************************************
\section{Perturbative ... or not?}
A quantum theory is modeled as a Feynman Process,
which is a representation of a Feynman Category (the causal structure),
defined bellow.

As explained in \S\ref{S:physmotiv},
representing a causal structure ``a la Markov'' is not a perturbative approach
``per se'', but rather a natural resolution of the quantum system,
called in \cite{I-DWT} a {\em Quantum Dot Resolution}.
This is the {\em homological approach} to modeling quantum systems.

On the other hand it definitely looks perturbative,
when the class of ``paths'' of the PROP is derived 
{\em perturbatively}, by expanding as a Taylor-Feynman series \cite{I-pqft}
a Gaussian (matrix) integral using Wick's (Matrix) Theorem. %

What is new in this homological approach to ``space-time'',
which {\em is} in fact a model for a discrete (quantum) ``space-time''
(like spin networks and foams of Loop Quantum Gravity etc.), 
is that it is NOT just an approximation of a continuum model; 
rather the continuum limit is the classical approximation of the 
quantum model.

we will briefly mention some additional features of
the homological approach to quantum modeling:

\quad 1)Enables the {\em resolution} / multi-scale analysis (MRA);

\quad 2) Enables in a fundamental way ``symmetry fluctuations''
of the causal structure and the natural implementation of
the concepts of entropy, information and the possibility
to involve them in mechanisms of mass generation;

\quad 3) It prompts for an {\em Q-information flow} interpretation
of the quantum dynamics (FPI at two levels: partition function AND
within a possible Feynman graph/RS \cite{I-QID}).

So, what is a {\em Feynman Category}?
Suggested by the DG-coalgebra structure of Feynman graphs \cite{I-pqft}
$$\xymatrix @C=3pc @R=1pc {
& [k] \dto \rto & \gamma \ar@{=>}[d] \rto & [l]\dto \\
\gamma\subset\Gamma\to \Gamma/\gamma 
& [n] \rto \dto & \Gamma \rto \ar@{=>}[d] & [m]\dto \\
& [n-k] \rto & \Gamma/\gamma \rto & [m-l],
}$$
a {\em Feynman Category} (FC) is a 
{\bf DG-coalgebra PROP} of finite type over $(N,+)$.

The coalgebra structure encodes the factorization of 
morphisms/processes in the {\em space-like/resolution depth} direction.

A {\em Feynman Process} is a representation of a Feynman Category. 
Feynman rules, Kontsevich's rules etc. are recipes to build such representations.
To have a physically meaningful action,
an additional $SU(2)$ or conformal action should be taken into account.

% **********************************************************
% 			page 15
% **********************************************************
\section{Homology and cohomology of Feynman Categories}
To study a causal structure one needs to build and study
Feynman Processes, i.e. the corresponding representation theory.

As mentioned before, a machinery able to generate Feynman Categories
is Wick's Theorem applied to Gaussian integrals and matrix models.

Once a Feynman Category $\C{F}$ is given, 
a general (classical) method for constructing Feynman Processes
is the so called {\em sigma model}.
An ambient (semi-)Riemannian space-time is used to embed graphs, 
Riemann surfaces, cobordisms,
and in general the ``paths'' of the geometric category,
in order to right an action functional (functorial):
$$``Hom(geometric\ category, ambient\ manifold)''$$
One can then study the {\bf generalized (co)homology of a manifold}:
the homology/homotopy of the corresponding configuration space.

QFT, CFT, TQFT etc. are, as Jim Stasheff put it so suggestively,
cohomological physics. These quantum theories are (classes of functors)
are generalized cohomology theories and the corresponding cohomologies
are essentially their classical limits (tree-level restrictions)
Precise statements may be found in \cite{Vor}.
For example, a (1+1)-TQFT (Frobenius algebra) at tree-level is a
commutative algebra; the cohomology of a string background is a (1+1)-TQFT 
(Frobenius algebra).

% **********************************************************
% 			From continuum to discrete
% **********************************************************
\section{From {\em Continuum} to {\em Discrete}}
But the key issue here is to develop an {\bf abstract theory} of FC,
by moving away from the {\em continuum}/manifold theory to the 
{\em discrete} world of (finite type) resolutions.

\subsection{Discrete manifolds}
To benefit from the rich theory already developed,
start from a manifold, and:

\vspace{.1in}
\quad 1) Discretize the manifold; 

\vspace{.1in}
\quad 2) {\bf Pullback ``The Theory''} to the ``geometric category''
(e.g. history of abstract groups: groups of transformations to abstract groups);

\vspace{.1in}
\quad 3) Study its representation theory 
(Generalized cohomology with coefficients in a modular category).

As a warm-up exercise, start with Feynman-Kontsevich graphs, 
and then ... attack String Theory!

The ``pull-back philosophy'' is a growing trend in high-energy physics:

\quad - Loop Quantum Gravity (LQG) starts from the general relativity 
manifold picture and ends up with a 
discrete space-time (spin-networks and spin-foams etc.);

\vspace{.1in}
\quad - String Theory as a ``background-free'' theory found an 
algebraic formulation (string backgrounds), yet a formulation
still regarding the representation theory;
an intrinsic formulation is a future project \cite{I-QID}.

What is by now clear, is that 
it saves time to ``adopt'' the ``Feynman Picture'' from the beginning. 
Feynman processes are ``just'' enriched and complexified Markov processes ...

\subsection{What is ``Space-Time''?}
Paraphrasing \cite{GM}:
``It does not matter; all we need is a resolution''!

The Kontsevich's sigma-model quasi-isomorphism {\em is}
such a candidate, at least at the level of type of construction,
of a {\em resolution} of the ``ambient space-time''.
More precisely, 
the Hochschild DGLA of the Poisson algebra of observables 
of a Poisson manifold $M$ is a {\em formal manifold},
representing an algebraic substitute for the original manifold
the same way in algebraic geometry 
the commutative algebra of functions determines the space.
Now the formality morphism plays the role of a universal resolution
of the formal manifold.

The study of a quantum theory from this point of view will be the 
subject of a different article.
But the general philosophy that a quantum theory is a Feynman Process
representing a causal structure modeled as a Feynman category
which is a ``resolution of space-time'' (external degrees of freedom
not ``fixed'', due to quantum fluctuations,
which in turn corresponds to our multi-scale analysis of the system),
is clear.
The rest is the mathematical study of the $\sigma$-model $Hom(\cdot, M)$,
as a {\em Configuration Functor} and its derived functors,
the {\bf FPI-quantization}.

In our opinion, 
Kontsevich's approach to deformation quantization based on graphs
is more important through its ideology and 
underlying algebraic structures present.

% *************************************************************
% 		Missing Physics
% *************************************************************
\section{``Missing'' Physics?}
% see pages 11,12 from the Bucharest talk
The computer science perspective on quantum physics hints
towards the equivalence between energy and information, at a fundamental level,
not just as a balance equation between energy and entropy as in thermodynamics.
Here ``energy'' includes matter, with its space-time determination via Einstein's 
principles of general relativity, and therefore this new,
yet to be precisely formulated, equivalence principle extends 
Einstein's equivalence principles $E=mc^2$ and $m_a=m_g$.

\subsection{Entropy is a measure of symmetry!}
The road to such a unification goes through 
Shannon entropy:
$$S=k \ln W\quad <=> \quad S=k |Aut(\Gamma)|$$
with important conceptual implications: {\em entropy is a measure of symmetry}! 

Now in Feynman Processes the causal structure fluctuates with the (homological) resolution degree,
and this includes its symmetry. 
The mathematical framework is adequate to implement these ideas as 
{\em Quantum Information Dynamics} (QID) \cite{I-QID}.

The straightforward way to implement the equivalence between mass-energy and entropy (quantum information),
at the fundamental level is by including entropy in the action, via the symmetry group:
$$Z=\int_{\gamma\in Hom(In,Out)} e^{i \C{S}(\gamma)}/|Aut(\gamma)|,$$
$$e^{i \C{S}(\gamma)}/|Aut(\gamma)|= e^{-\ln |Aut(\gamma)|+i \C{S}(\gamma)}$$
$$\framebox{$Z=\int\limits_{Hom(In,Out)} e^{S+i\C{S}}d\mu$}.$$
This has already the flavor of Green functions $S+i\C{S}=-\ln |Aut(\gamma)|+i \C{S}(\gamma)$,
towards the interpretation of Riemann surfaces as ``networks'' propagating quantum information
(amplitudes/CFT). 
It is an avenue prone for a true merger of Euclidean formulation and the statistical formulation of
QFT, by generalizing {\bf Wick rotation} (space-like versus time-like description),
to achieve a complete {\bf unification of space and time} \cite{I-QID}.

The fluctuation of the symmetry includes the classical breaking of symmetry (change in entropy),
with its mathematical formulation as a mass generation mechanism.
Ideologically speaking, mass and gravity are entropic effects.

\subsection{Additional ``evidence''}
That there is such a level of unification between 
energy/matter/space-time and quantum information
is already apparent in the quantum radiation laws of black holes.
One of these fundamental laws stipulate that entropy is proportional 
with the surface area \cite{I-UP}.
As Lee Smolin puts it, this law suggests that the surface is discreet (quantum)
and the unit is a ``pixel of space-time''.
Now since space-time is matter (quantum fields'' etc.),
and since the unit of quantum information is the qubit (spin'' is less striking!),
we arrive at an equivalence relation; yet the symmetry contribution is still missing ...

There are other connections, or rather supporting evidence towards such a unification,
as well as for the conceptual benefits (power of a comprehensive model).
We will only mention the concept of quantum potential in B. J. Hiley's work,
or better (Green's) quantum information potential,
and Bob Coecke's quantum {\bf information-flow} at the level of QM, 
meaning in our view, the quantum computation ``order'' (flow),
tightly related with the space-time coordinatisation, 
as we will explain it in \cite{I-QID}.

The author's personal feeling is that there is ``new physics'' at the horizon:
the rise of ``Low Entropy Physics'' (LEP) versus 
the down of ``High-Energy Physics'' (HEP-th) \cite{I-DWT}.
The mathematics is already present ... and the interpretation, finally emerges!

% **********************************************************
% 			Conclusions
% **********************************************************
\section{Conclusions}
HEP-th is a study of {\bf representations of PROPs}:
Feynman, Segal, cobordism categories etc.

\subsection{``New Mathematics''?}
No, not really!
A {\bf Feynman Category} (PROP),
i.e. a ``geometric category'' playing the role of causal structure (e.g. Feynman graphs),
is a homological resolution of whatever ``space-time'' is.
It usually ``looks perturbative'', being obtained via Wick's Theorem 
from Gaussian integrals (matric models),
but conceptually is a homological resolution enabling multi-scale analysis.

Then the tree-level of QFT, CFT, String Theory and other algebras over Feynman PROP
can be thought of as ``derived functors'' of a {\bf Configuration Functor} 
in the context of a $\sigma$ model.
From this point of view Kontsevich's construction appears more important
through the interpretation of the formality quasi-isomorphism as
an augmentation of a resolution of the sigma model,
representing the underlying causal structure:
$$Feynman\ Category \quad \overset{quasi-iso}{\to}\quad Formal Manifold$$
$$DGCA/DGLA:\quad (\C{G},d,\Delta) \overset{\epsilon}{\to} (End(T(A)),Q).$$

\subsection{``New Physics''?}
May be!
Entropy as a measure of symmetry enters the dynamics of the
``space-time'' (causal structure) in the context of 
quantum information dynamics in a extended sense:
$$QID: \quad Feynman\ Processes\ as\ Quantum\ Information\ Dynamics.$$ 
Gravity (quantum, as everything else), 
is ``just'' an organizing principle naturally emerging from the 
relevant (grand) unification.

%=======================================================================
%                      Bibliography
%=======================================================================
\bibliographystyle{my-h-elsevier}

\end{document}